\documentclass{jgcc} 

\pdfoutput=1
\usepackage{lastpage}
\jgccdoi{15}{1}{3}{11315}  
\jgccheading{}{\pageref{LastPage}}{}{}{May~15,~2023}{Oct.~19,~2023}{}    


\newcommand{\bi}{\bibitem}
\newcommand{\nb}{\newblock}

\newcommand{\be}[1]{\begin{equation}\label{#1}}
	\newcommand{\ee}{\end{equation}}

\newcounter{ppp}
\newcommand{\la}{\langle\,}
\newcommand{\ra}{\,\rangle}

\newcommand{\bdelta}{\bar\delta}
\newcommand{\distan}{\mathop{\rm dist}}

\newcommand{\hgt}{\mathop{\rm ht}}

\usepackage{hyperref}

\begin{document}

\title[Amenability problem for Thompson's group $F$]{Amenability problem for Thompson's group $F$: \\state of the art}

\author[V.~Guba]{Victor Guba}	
\address{Vologda State University,
	15 Lenin Street,
	Vologda,
	Russia
	160600}	
\email{gubavs@vogu35.ru}  
\thanks{This work is supported by the Russian Science Foundation, project no. 23-21-00289.}	




\begin{abstract}
  \noindent 	This is a survey of our recent results on the amenability problem for Thompson's group $F$. They mostly concern esimating the density of finite subgraphs in Cayley graphs of $F$ for various systems of generators, and also equations in the group ring of $F$. We also discuss possible approaches to solve the problem in both directions.
\end{abstract}

\maketitle

\section{Preliminaries}
\label{prel}

In this survey we collect a number of our recent results on the amenability problem for Thompson's group $F$. They are mostly contained in~\cite{Gu21a,Gu21b,Gu21c,Gu22s,Gu22se,Gu22,Gu23,Gu23b}. In the this section we introduce concepts and notation used throughout the paper. The reader who is familiar with that can skip this background.

We recall the concept of a graph in the sense of Serre, Cayley graph of a group (right anf left), internal and external vertices of a subgraph of a Cayley graph, and various kinds of its boundary. It is convenient to use automata terminology saying that an automaton accepts (or does not accept) a generator. We introduce an important concept of density of Cayley graphs and their subgraphs together with the related concept of their (Cheeger) isoperimetric constant.

Then we recall the definition of amenabilty of semigroups and groups, discuss Folner criterion and some other necessary and sufficient conditions for amenabilty. In order to work with elements of Thompson's group $F$ whose basic properties are discussed in a separate subsection, we present terminology related to rooted binary trees and forests. We recall the definition of the group $F$ and list some of its important properties. The amenability problem for this group is in the centre of our attention.

Section~\ref{dencg} contains the description of our recent result on the density of finite subgraphs of the Cayley graph of $F$ in various generating systems. One of the most important results is Theorem~\ref{over3.5} where we show that the density of the Cayley graph of $F$ in standard generators strictly exceeds $3{.}5$. This is an improvement of the Belk-Brown construction for which there was a conjecture of its optimality.

Section~\ref{eqgr} concerns equations in the group ring of $F$ and their systems. Due to results by Bartholdi and Kielak, we know that the amenabilty problem for $F$ is equivalent to the Ore condition for group rings, that is, the question whether equation $au=bv$ in this ring has nonzero solutions. We present a number of results in this direction giving a reduction of the problem to the case of monoid rings $K[M]$ of the positive monoid $M\subset F$ and homogeneous polynomials $a$, $b$. We solve some partial cases of this problem including the case of linear polynomials. We also give a description for the set of solutions of some important equations and their systems.

Section~\ref{misc} contains a new description of exhausting finite subsets of the Cayley graph of $F$ and a new algorithm to solve the word problem in this group. Besides, we discuss possible approaches to the amenability problem to solve it in both directions.

\subsection{Cayley graphs}
\label{cay}

By a {\em graph} we mean a (non-oriented) graph in the sense of Serre~\cite{Se80}. This means that each geometric edge consists of two mutually inverse directed edges. As a formal concept, this is a 5-tuple $\Gamma=\la V,E,^{-1},\iota,\tau\ra$, where $V$, $E$ are disjoint sets (of {\em vertices} and {\em edges}, respectively); $^{-1}\colon E\to E$, $\iota\colon E\to V$, $\tau\colon E\to V$. We assume that $e\ne e^{-1}$, $(e^{-1})^{-1}=e$, $\iota(e^{-1})=\tau(e)$, $\tau(e^{-1})=\iota(e)$ for all $e\in E$. Here $e^{-1}$ is called the {\em inverse} edge of $e$, $\iota(e)$ is the {\em initial} vertex of $e$, $\tau(e)$ is the {\em terminal} vertex of $e$.


Let $A$ be an {\em alphabet}, that is, a nonempty set of {\em symbols} (or {\em letters}). Taking a disjoint bijective copy $A^{-1}$ of the set $A$, we get the set $A^{\pm1}=A\cup A^{-1}$ called a {\em group alphabet}. It has an involution $^{-1}$ without fixed points defined by $(a^{-1})^{-1}\rightleftharpoons a$. Elements of the free monoid $(A^{\pm1})^{\ast}$ are called {\em group words}.

Let $G$ be a group equipped by mapping $A\to G$ such that the image of $A$ generates $G$. We will say that $A$ is a set of (group) generators for $G$. Usually $A$ is identified with its image provided the above map is injective. For more general needs it is convenient to include the case when different symbols may denote the same group generator. In this case we may say that $A$ defines a multiset of generators, say, $\{a_1=x_0,a_2=x_0,a_3=x_1\}$, where $x_0$ is doubly repeated.

Let $G$ be a group generated by $A$ in the above sense. Its {\em right Cayley graph} $\Gamma_r={\mathcal C}(G;A)$ is defined as follows. The set of vertices is $G$; for any $a\in A^{\pm1}$ we put a directed edge $e=(g,a)$ with $\iota(e)=g$, $\tau(e)=ga$. Its inverse is $e^{-1}=(ga,a^{-1})$. Each edge has a {\em label}; for the above edge $e=(g,a)$ it is just $a$. The label function can be naturally continued to the set of paths. 


The concept of a {\em left Cayley graph} $\Gamma_l={\mathcal C}(G;A)$ is defined in a similar way. Here the set of vertices is also $G$. A directed egde $e=(a,g)$ with label $a$ has $\iota(e)=ag$, $\tau(e)=g$. That is, going along the edge labelled by $a$, means cancelling $a$ on the left. The inverse egde is $e^{-1}=(a^{-1},ag)$. Labels of paths are defined similarly for this case.

\subsection{Automata, density, and isoperimetric constants}
\label{aut}

The cardinality of a finite set $Y$ will be denoted by $|Y|$. 

Let $G$ be an infinite group generated by $A$. For our needs we assume that $A$ is always finite, $|A|=m$.  Let $\Gamma={\mathcal C}(G,A)$ be the Cayley graph of $G$, right or left. To any finite nonempty subset $Y\subset G$ we assign a subgraph in $\Gamma$ adding all edges connecting vertices of $Y$. So given a set $Y$, we will usually mean the corresponding subgraph. This is a labelled graph that we often call an {\em automaton}. For each $g\in Y$ we have exactly $2m$ directed edges in $\Gamma$ starting at $g$, where $a\in A^{\pm1}$. If the endpoint of such an edge with label $a$ belongs to $Y$, then we say that the vertex $g$ of our automaton $Y$ {\em accepts} $a$. For the case of right Cayley graphs this means $ga\in Y$, for the case of left Cayley graphs this means $a^{-1}g\in Y$. 

A vertex $g\in Y$ is called {\em internal} whenever it accepts all labels $a\in A^{\pm1}$. That is, the degree of $g$ in $Y$ equals $2m$. Otherwise we say that $g$ belongs to the {\em inner boundary} of $Y$ denoted by $\partial Y$. 

By $\distan(u,v)$ we denote the distance between two vertices in $\Gamma$, that is, the length of a shortest
path in $\Gamma$ that connects vertices $u$, $v$. For any vertex $v$ and a number $r\ge0$ let $B_r(v)$ denote the ball of radius $r$ around $v$, that is, the set of all vertices in $\Gamma$ at distance $\le r$ from $v$. For any set $Y$ of vertices, by $B_r(Y)$ we denote the $r$-neighbourhood of $Y$, that is, the union of all balls $B_r(v)$, where $v$ runs over $Y$. By
$\partial_o Y$ we denote the {\em outer boundary\/} of $Y$, that is, the set $B_1(Y)\setminus Y$.

An edge $e$ is called {\em internal} whenever it connects two vertices of $Y$. If an edge $e$ connects a vertex of $Y$ with a vertex outside $Y$, then we call it {\em external}. That is, $e$ connects a vertex in $\partial Y$ with a vertex in $\partial_o Y$. The set of external edges form the {\em Cheeger boundary} of $Y$ denoted by $\partial_{\ast}Y$. 

To any vertex $v$ in $\partial Y$ we can assign an external edge starting at $v$. This gives an injection from  $\partial Y$ to  $\partial_{\ast}Y$. On the other hand, there exist at most $2m$ edges in $\partial_{\ast}Y$ staring at $v$. This implies inequalities  $|\partial Y|\le|\partial_{\ast}Y|\le2m|\partial Y|$. The same arguments show that $|\partial_o Y|\le|\partial_{\ast}Y|\le2m|\partial_oY|$. 

By the {\em density} of a subgraph $Y$ we mean its average vertex degree. This concept was introduced in~\cite{Gu04}; see also~\cite{Gu21a}. It is denoted by $\delta(Y)$. A {\em Cheeger isoperimetric constant} of the subgraph $Y$ is the quotient $\iota_*(Y)=|\partial_{\ast}Y|/|Y|$. It follows directly from the definitions that $\delta(Y)+\iota_*(Y)=2m$. Indeed, each vertex $v$ has degree $2m$ in the Cayley graph $\Gamma$. This is the sum of the number of internal edges starting at $v$, which is $\deg_Y(v)$, and the number of external edges starting at $v$. Taking the sum over all $v\in Y$, we have $2m|Y|$, which is equal to $\sum\limits_v\deg_Y(v)+|\partial_{\ast}Y|$. Dividing by $|Y|$, we get the above equality.

(We can mention also the paper~\cite{AGL08} where some computaional experiments with densities of subgraphs are presented.)

By the density of the Cayley graph $\Gamma={\mathcal C}(G;A)$ we mean the number $\bdelta(\Gamma)=\sup\limits_Y\delta(Y)$, where $Y$ runs over all nonempty finite subgraphs. Analogously, the Cheeger isoperimetric constant of the group $G$ with respect to generating set $A$, is the number $\iota_*(\Gamma)=\inf\limits_Y\iota_*(Y)$. Clearly, $\bdelta(\Gamma)+\iota_*(\Gamma)=2m$.

In further sections, working with isoperimetric constants, we will write $\iota_*(G;A)$ instead of the above notation.

\subsection{Amenability of semigroups and groups}
\label{amen}

Let $S$ be a semigroup. Suppose that there exists a mapping $\mu\colon{\mathcal P}(S)\to[0,1]$ from the power set of $S$ into the unit interval satisfying the following conditions:
\vspace{0.5ex}

1) $\mu$ is additive, that is, $\mu(A\cup B)=\mu(A)+\mu(B)$ for any disjoint subsets $A,B\subseteq S$;

2) $\mu$ is left invariant, that is, $\mu(sA)=\mu(A)$ for any $s\in S$, $A\subseteq S$;

3) $\mu$ is normalized, that is, $\mu(S)=1$.
\vspace{0.5ex}

Then the semigroup $S$ is called {\em left amenable}.
\vspace{1ex}

The definition of {\em right amenable} semigroups is given in a similar way. These concepts differ in general. However, for the case of groups they are equivalent. Moreover, both of them will be equivalent to the conditions with the two-sided invariance. The proof can be found in \cite{GrL}.

We do not list all well-known properties of (non)amenable groups. It is sufficient to refer to one of modern surveys like~\cite{Sap14}. Just notice that all finite and abelian groups are amenable. The class of amenable groups is closed under taking subgroups, homomorphic images, group extensions, and directed unions of groups. The groups in the closure of the classes of finite and abelian groups under this list of operations are called {\em elemetary amenable} (EA). Also we need to say that free groups of rank $> 1$ are not amenable. 
\vspace{1ex}

We will often refer to the following F\o{}lner criterion~\cite{Fol}. Here we restrict ourselves to the case of finitely generated groups.

\begin{prop}
	\label{fol}
	A group $G$ with finite set of generators $A$ is amenable if and only if its Cheeger isoperimetric constant iz zero: $\iota_*(G;A)=0$.
\end{prop}

This holds for any finite set of generators. Equivalently, one can say that the density of the Cayley graph $\Gamma={\mathcal C}(G;A)$ (right or left) has its maximum value: $\bdelta(\Gamma)=2m$, where $m=|A|$.

In practice, to establish amenability of a group, it is sufficient to construct a collection of nonempty finite subgraphs $Y$ in the Cayley graph such that $\inf\limits_Y|\partial_{\ast}Y|/|Y|=0$. Such subsets of vertices are called {\em Folner sets}. Informally, this means that almost all vertices of these sets are internal.

If one expects that a group is not amenable, or the answer is unknown, it is reasonable to construct sets with the least possible value of $\iota_*(G;A)$. We will discuss this approach later.

In case when we are going to establish non-amenability of a group, it is useful to apply some other criteria.

Let $G$ be a group generated by a finite set $A$. A {\em doubling function} on $G$ is a mapping $\psi\colon G\to G$ such that

a$)$ for all $g\in G$ the distance $\distan(g,\phi(g))$ is bounded from above by a constant $K > 0$,

b$)$ any element $g\in G$ has at least two preimages under $\psi$.

\begin{prop}
	\label{grom}
	A group $G$ with finite set of generators $A$ is non-amenable if and only if it admits a doubling function.
\end{prop}

This criterion is often attributed to Gromov. An elegant proof of it can be found in \cite{CGH}, see also \cite{DeSS}. Note that this property also does not depend on the choice of a finite generating set. 

An interesting partial case happens if the constant $K$ in the above definition equals $1$. In~\cite{Gu04} we proved that for a 2-generator group $G$, a doubling function with constant $K=1$ exists if and only if the density of the Cayley graph does not exceed $3$. Equivalently, one can say that the Cheeger isoperimetric constant is at least one: $\iota_*(G;A)\ge1$. 

We call a group $G$ {\em strongly non-amenable} (with respect to a given finite generating set $A$) whenever this inequality holds: $\iota_*(G;A)\ge1$. It is interesting to find a nesessary and sufficient condition for this property for groups with any finite number of generators. This will be done in Section~\ref{dencg}.
\vspace{1ex}

Another practical criterion for non-amenabilty can be stated in terms of flows on Cayley graphs. Let us introduce the terminology for that.

A {\em flow} on a graph $\Gamma$ is a real-valued function $f\colon E\to\mathbb R$ such that $f(e^{-1})=-f(e)$. We say that $f(e)$ is the flow through the edge $e$. Given a vertex $v$, we define an {\em inflow} to it as a sum of flows through all edges with $v$ as a terminate vertex.

The following criterion in terms of flows is essentially known. It can be derived from the above criterion in terms of doubling functions.

\begin{prop}
	\label{critfl}
	Let $G$ be a group with finite generating set $A$, and let $\Gamma={\mathcal C}(G;A)$ be its Cayley graph. The group $G$ is non-amenable if and only if there exist constants $C > 0$ and $\varepsilon > 0$, and a flow $f$ on $\Gamma$ with the following properties:
	
	a$)$ The absolute value of the flow through each edge is bounded: $|f(e)|\le C$ for all $e\in E$;
	
	b$)$ The inflow into each vertex is at least $\varepsilon$.	
\end{prop}

\subsection{Rooted binary trees and forests}
\label{rbtf}

We add this short subsection to introduce some notation used in the paper. 

Formally, a {\em rooted binary tree} can be defined by induction.

1) A dot $\cdot$ is a rooted binary tree.

2) If $T_1$, $T_2$ are rooted binary trees, then $(T_1\hat{\ \ }T_2)$ is a rooted binary tree.

3) All rooted binary trees are constructed by the above rules.
\vspace{1ex}

Instead of formal expressions, we will use their geometric realizations. A dot will be regarded as a point. It coincides with the root of that tree. If $T=(T_1\hat{\ \ }T_2)$, then we draw a {\em caret\/} for $\hat{}$ as a union of two closed intervals $AB$ (goes left down) and $AC$ (goes right down). The point $A$ is the {\em root} of $T$. After that, we draw trees for $T_1$, $T_2$ and attach their roots to $B$, $C$ respectively in such a way that they have no intersection. It is standard that for any $n\ge0$, the number of rooted binary trees with $n$ carets is equal to the $n$th Catalan number $c_n=\frac{(2n)!}{n!(n+1)!}$.

Each rooted binary tree has {\em leaves\/}. Formally, they are defined as follows: for the one-vertex tree
(which is called {\em trivial\/}), the only leaf coincides with the root. In case $T=(T_1\hat{\ \ }T_2)$, the set of leaves
equals the union of the sets of leaves for $T_1$ and $T_2$. In this case the leaves are exactly vertices of degree
$1$.

We will also need the concept of a {\em height\/} of a rooted binary tree. For the trivial tree, its height equals
$0$. For $T=(T_1\hat{\ \ }T_2)$, its height is $\hgt T=\max(\hgt T_1,\hgt T_2)+1$.

Now we define a {\em rooted binary forest\/} as a finite sequence of rooted binary trees $T_1$, ... , $T_m$,
where $m\ge1$. The leaves of it are the leaves of the trees. It is standard from combinatorics that the number
of rooted binary forests with $n$ leaves also equals $c_n$. The trees are enumerated from left to right and they
are drawn in the same way.

A {\em marked\/} (rooted binary) forest is a (rooted binary) forest where one of the trees is distinguished.

\subsection{Thompson's group $F$}
\label{thgf}

We define the Richard Thompson group $F$ in a combinatorial way, using the following infinite
group presentation

\be{xinf}
\la x_0,x_1,x_2,\ldots\mid x_j{x_i}=x_ix_{j+1}\ (i<j)\,\ra.
\ee
This group was found by Richard J. Thompson in the 60s. We refer to the survey \cite{CFP} for details. (See also \cite{BS,Bro,BG}.) A recent survey on the subject with respect to the amenability problem of $F$ can be found in~\cite{Gu22s}. 

It is easy to see from the relations of~(\ref{xinf}) that for any $n\ge2$, one has $x_n=x_0^{-(n-1)}x_1x_0^{n-1}$ so
the group is generated by $x_0$, $x_1$. It can be given by the following presentation with two defining relations

\be{x0-1}
\la x_0,x_1\mid x_1^{x_0^2}=x_1^{x_0x_1},x_1^{x_0^3}=x_1^{x_0^2x_1}\ra,
\ee
where $a^b=b^{-1}ab$ by definition. Also we define a commutator $[a,b]=a^{-1}a^b=a^{-1}b^{-1}ab$
and notation $a\leftrightarrow b$ whenever $a$ commutes with $b$, that is, $ab=ba$.

Each element of $F$ can be uniquely represented by the {\em normal form\/}, that is, an expression of the form
\be{nf}
x_{i_1}x_{i_2}\cdots x_{i_s}x_{j_t}^{-1}\cdots x_{j_2}^{-1}x_{j_1}^{-1},
\ee
where $s,t\ge0$, $0\le i_1\le i_2\le\cdots\le i_s$, $0\le j_1\le j_2\le\cdots\le j_t$ and the following is true: if (\ref{nf}) contains both $x_i$ and $x_i^{-1}$ for some $i\ge0$, then it also contains $x_{i+1}$ or $x_{i+1}^{-1}$ (in particular, $i_s\ne j_t$).
\vspace{1ex}

Equivalent definitions of $F$ can be given in terms of piecewise-linear functions. Although these definitions are very popular, we do not describe it here since we will not use them in our survey. 

It is known from~\cite{GbS} that $F$ is the diagram group over the simplest semigroup presentation$\la x\mid x=x^2\ra$. This way to represent elements of $F$ is very useful for many situations. Sometimes it is preferable to use non-spherical diagrams over the same presentation instead of spherical ones. The latter approach was described in~\cite{Gu04}. Diagrams can be replaced by dual graphs. This leads to a standard way to represent elements of $F$ as pairs of rooted binary trees. For the needs of this survey, it suffices to work with elements of the positive monoid $M$ repersenting them as marked binary forests according to Section~\ref{rbtf}.
\vspace{1ex}

The group $F$ has no free non-abelian subgroups~\cite{BS}. It is known~\cite{Chou} that $F$ is not elementary amenable. However, the famous problem about amenability of $F$ remains open. The question whether $F$ is amenable was asked by Ross Geoghegan in 1979; see~\cite{Geo,Ger87}. There is no common opinion on the answer: some of specialists in this area try to prove non-amenability, some of them believe that the group is amenable. There is a number of papers with attempts to solve the problem in both directions. The author always believed in non-amenability of $F$ trying to prove this property. Now we are not sure on that answer, the reasons will be explained later.

If $F$ is amenable, then it is an example of a finitely presented amenable group, which is not EA. If it is not
amenable, then this gives an example of a finitely presented group, which is not amenable and has no free non-abelian subgroups.  Note that the first example of a non-amenable group without free non-abelian subgroups has been constructed by Ol'shanskii~\cite{Olsh}. The question about such groups was formulated in \cite{Day}, it is also often attributed to von Neumann~\cite{vNeu}. Adian~\cite{Ad83} proved that free Burnside groups with $m>1$ generators of odd exponent
$n\ge665$ are not amenable. The first example of a finitely presented non-amenable group without free non-abelian subgroups has been constructed by Ol'shanskii and Sapir~\cite{OlSa}. Grigorchuk~\cite{Gri} constructed the first example of a finitely presented amenable group not in EA.
\vspace{2ex}

It is not hard to see that $F$ has an automorphism given by $x_0\mapsto x_0^{-1}$, $x_1\mapsto x_1x_0^{-1}$. To check that, one needs to show that both defining relators of $F$ in (\ref{x0-1}) map to the identity. This is an easy calculation using normal forms. After that, we have an endomorphism of $F$. Aplying it once more, we have the identity map. So this is an (outer) automorphism of order $2$.


Later we will add more arguments to the importance of the symmetric set $S=\{x_1,\bar{x}_1\}$, where $\bar{x}_1=x_1x_0^{-1}$. Obviously, $S$ also generates $F$. We let $\alpha=x_1^{-1}$, $\beta=\bar{x}_1^{-1}=x_0x_1^{-1}$.  Applying Tietze transormations to~(\ref{x0-1}) one can get a presentation of $F$ in the new generating set:

\be{albet}
\la\alpha,\beta\mid\alpha^{\beta}\leftrightarrow\beta^{\alpha},\alpha^{\beta}\leftrightarrow\beta^{\alpha^2}\ra.
\ee

From the symmetry reasons we know that $\beta^{\alpha}\leftrightarrow\alpha^{\beta^2}$ also holds in $F$. Therefore, it is a consequence of the two relations of~(\ref{albet}). Moreover, one can check that for any positive integers $m$, $n$ it holds $\alpha^{\beta^m}\leftrightarrow\beta^{\alpha^n}$ as a consequence of the defining relations.
\vspace{1ex}

We will often work with a positive monoid $M$ of the group $F$. It is defined by the monoid presentation that coincides with~(\ref{xinf}). The group $F$ is the group of quotients of this monoid so that $F=MM^{-1}$. Elementary reasons show that any finite subset in $F$ can be moved into $M$ up to a right multiplication by an element $g\in F$ (see Lemma~\ref{gig} in Section~\ref{eqgr}). 

\section{Density of Cayley graphs}
\label{dencg}

In~\cite{Gu04} we constructed a family of finite subgraphs of the Cayley graph of $F$ in the standard generating set $\{x_0,x_1\}$. The densities of these subgraphs approach $3$. In the Addendum to the same paper we demonstrated a modification showing that the densities of finite subgraphs can strictly exceed $3$.

An essential improvement was made by Belk and Brown\cite{Be04,BB05}. They gave a family of finite subgraphs depending on two integer parameters $k\ge0$ and $n\ge1$. Here the density of the corresponding subgraphs in the same Cayley graph approaches $3{.}5$.

We need an explicit definition of these sets. Let $BB(n,k)$ denote the set of all marked binary forests with $n$ leaves, where all trees have height at most $k$. We regard it as a subset of the left Cayley graph of $F$ in standard generators. Let us describe how the group generators act on the vertices (all actions are left partial ones).

The generator $x_0$ acts by shifting the marker one position left if this is possible. Action of $x_0^{-1}$ means moving the marker one step to the right. The action of $x_1$ is as follows. If the marked tree is trivial, this is not applied. If the marked tree is $T=(T_1\hat{\ \ }T_2)$, then we remove its caret and mark the tree $T_1$. It is easy to see that applying $\bar{x}_1=x_1x_0^{-1}$ means the same replacing $T_1$ by $T_2$ for the marked tree (notice that $x_1$ acts first).

Now one can see that $x_1$ and $\bar{x}_1$ are totally symmetric. They generate $F$ so one can regard them as
the most natural generators besides the standard ones.

The action of $x_1^{-1}$ and $\bar{x}_1^{-1}$ is defined analogously. Namely, if the marked tree of a forest is
rightmost, then $x_1^{-1}$ cannot be applied. Otherwise, if the marked tree $T$ has a tree $T''$ to the right of
it, then we add a caret to these trees and the tree $(T\hat{\ \ }T'')$ will be marked in the result. Notice that
if we are inside $BB(n,k)$, then both trees $T$, $T''$ must have height $< k$: otherwise $x_1^{-1}$ cannot be
applied. As for the action of $\bar{x}_1^{-1}$, it cannot be applied if $T$ is leftmost. Otherwise the marked tree $T$ has
a tree $T'$ to the left of it. Here we add a caret to these trees and the tree $(T'\hat{\ \ }T)$ will be marked in the
result. As above, both trees $T'$, $T$ must have height $< k$ to be possible to stay inside $BB(n,k)$.
\vspace{1ex}

Let us define a sequence of polynomials by induction:
\be{phi0}
\Phi_0(x)=x
\ee
\be{phik}
\Phi_k(x)=x+\Phi_{k-1}(x)^2\quad\quad k\ge1.
\ee

Notice that $\Phi_k(x)$ is the generating polynomial for the set of trees of height at most $k$. This means that the coefficient on $x^n$ in this polynomial shows the number of such trees with $n$ leaves. This follows directly from~(\ref{phik}). The summand $x$ corresponds to the trivial tree (with one leaf); for the tree $T=(T_1\hat{\ \ }T_2)$ of height $\le k$ we have height $\le k-1$ for each of the trees $T_1$, $T_2$. By induction, the pair of them has generating function $\Phi_{k-1}(x)^2$. This agrees with~(\ref{phik}).
\vspace{1ex}

It is easy to see that the equation $\Phi_k(x)=1$ has a unique positive root that we denote by $\xi_k$. It is known that $\xi_k$ approaches $\frac14$ as $k\to\infty$. It is not hard to see that for a random marked binary forest from $BB(n,k)$, the probability to accept $x_0$ or $x_0^{-1}$ approaches $1$ as $n\to\infty$. As for $x_1$ and $x_1^{-1}$, the probability for a random vertex in the automaton $BB(n,k)$ to accept it approaches $\frac14$ for $n\gg k\gg1$. The same holds for symmteric generators $\bar{x}_1$ and $\bar{x}_1^{-1}$. 

It follows from these remarks that the density of the Cayley graph of $F$ in standard generators $\{x_0,x_1\}$ is at least $3{.}5$ since $x_0^{\pm1}$ are almost always accepted, and each of the $x_1^{\pm1}$ is accepted with probability close to $\frac14$. In other terms, it follows that $\iota_*(F,A)\le\frac12$ for $A=\{x_0,x_1\}$. This was a remarkable result obtained by Belk and Brown in~\cite{Be04,BB05}.

If we look from the same point of view to the symmetric generating set $\{x_1,\bar{x}_1\}$, then we see that the density of finite subsets $BB(n,k)$ will approach $3$ for $n\gg k\gg1$. However, this is not the best estimate.

Let we have a marked forest of the form $\ldots,T_{-1},T_0,T_1,\ldots$, where $T_0$ is marked. Suppose that $T_0$ is a trivial tree and each of the neighbour trees $T_{-1}$, $T_1$ has height $k$. In this case no generator of the form $x_1^{\pm1}$ or $\bar{x}_1^{\pm1}$ can be accepted by such a vertex in the automaton. This means that such vertices are isolated in $BB(n,k)$ as a subgraph of the left Cayley of $F$ in the symmetric generating set. The event to get an isolated vertex holds with guaranteed probability $p > 0$, where $p$ is a constant. Therefore, if we remove such vertices from the automaton, then the density will necessarily increase. Here is the main result from~\cite{Gu21a}.

\begin{thm}[\cite{Gu21a}]
\label{denssym}
The density of the Cayley graph of Thompson's group $F$ in symmetric generating set $S=\{x_1,\bar{x}_1\}$
is strictly greater than $3$.
\end{thm}

Notice that this trick give nothing for increasing density in case of standard generating set. Indeed, almost all vertices in $BB(n,k)$ accept both $x_0$ and $x_0^{-1}$. So if we remove a vertex isolated in the previous graph, we will need to remove 4 directed edges incident to it. This never has an effect of increasing density.

Let us mention one more result of the same paper. It concerns another symmetric set with the generator $x_0$ added. Notice that this fact was improved in our latest papers.

\begin{thm}
	\label{double}
	For the symmetric generating set $S=\{x_0,x_1,\bar{x}_1\}$, there exists finite subsets $Y\subset F$ in the
	Cayley graph of Thompson's group $F$ such that $|\partial Y| < |Y|$.
	
	Equivalently, the generating set $S$ does not have doubling property, that is, there are finite subsets $Y$ in
	$F$ such that the $1$-neighbourhood $\mathcal N_1(Y)=(S^{\pm1}\cup\{1\})Y$ has cardinality strictly less than $2|Y|$.
\end{thm}

In~\cite{Gu22} we introduced the concept of an evacuation scheme on a Cayley graph. Let we have an infinite group $G$ generated by a finite set $A$. Let $\Gamma={\mathcal C}(G;A)$ be its Cayley graph (right or left). To each vertex $v$ we assign an infinite simple path $p_v$  starting at $v$ in the Cayley graph. Suppose that there exists a constant $C$ such that each directed egde $e$ can participate in these paths at most $C$ times. In this case we say that the family $(p_v)_{v\in G}$ is an {\em evacuation scheme} on the Cayley graph $\Gamma$.

To be more precise: we claim that the total number of occurrences of each edge $e$ in paths of the form $p_v$ $(v\in G)$ does not exceed $C$.

Roughly speaking, the paths $p_v$ bring all verties to infinity. Without any restrictions, such an object always exists. However, if we claim that each edge participates in the scheme a uniformly bounded number of times, then we get the property equivalent to non-amenabilty. The following statement easily follows from known criteria. 

\begin{prop}[\cite{Gu22}]
	\label{critev}
	A group $G$ with finite generating set $A$ is non-amenable if and only if there exists an evacuation scheme on its Cayley graph.
\end{prop}

Suppose that we have an evacuation scheme with constant $K$ on the Cayley graph $\Gamma={\mathcal C}(G;A)$. Let $Y$ be a finite nonempty subset of $G$. We know that each path $p_v$ ($v\in Y$) must leave $Y$ at some step. Hence there exists an initial segment $\bar{p}_v$ of $p_v$ that has its terminal point on $\partial Y$. This leads to the concept of an {\em evacuation scheme with constant $K$ on a finite subgraph}. This is a collection of paths in $Y$ of the form $\bar{p}_v$. This finite path starts at $v$ and ends on the inner boundary $\partial Y$. For each edge $e$ we claim that the total number of its occurrences in the paths $\bar{p}_v$ ($v\in Y$) does not exceed $K$.

Having an evacuation scheme on $Y$, we can assume that paths $\bar{p}_v$ are simple (otherwise we can remove some loops). Also we can claim that if $e$ occurs in the evacuation paths, then $e^{-1}$ does not occur. Indeed, if $\bar{p}_v=p_1ep_2$, $\bar{p}_u=p_3e^{-1}p_4$, then one can replace these evacuation paths by $p_1p_4$, $p_3p_2$, respectively. 

The following fact easily follows from K\"onig's Lemma~\cite{Kon}.

\begin{prop}[\cite{Gu22}]
\label{kon}
Let $G$ be an infinite group generated by a finite set $A$. An evacuation scheme on the Cayley graph $\Gamma={\mathcal C}(G;A)$ exists if and only if there exist evacuation schemes on all its finite nonempty subgraphs. The constants for both cases are the same.
\end{prop}

An important case of an evacuation scheme happens if $C=1$ in the definition. In this case we say that we have a {\em pure evacuation scheme} on the Cayley graph. This means that each directed edge can participate in evacuation paths of the form $p_v$ at most once. We obtained a criterion for existing of such a scheme.

\begin{thm}[\cite{Gu22}]
	\label{main2}
	Let $G$ be an infinite group generated by a finite set $A$. A pure evacuation scheme on the Cayley graph $\Gamma={\mathcal C}(G;A)$ exists if and only if the group $G$ is strongly non-amenable with respect to $A$, that is, its Cheeger isoperimetric constant is at least one: $\iota_*(G;A)\ge1$.	
\end{thm}

The following statement improves one of our previous results.

\begin{thm}[\cite{Gu22}]
	\label{isopsym0}
	The Cheeger isoperimetric constant of the Cayley graph of Thompson's group $F$ in extended symmetric generating set $\{x_0,x_1,\bar{x}_1\}$ is strictly less than $1$.
\end{thm}

According to Theorem~\ref{main2}, this means that there is no pure evacuation scheme on the Cayley graph of $F$ in these generators.
\vspace{1ex} 

The construction used in the proof of Theorem~\ref{main2} was almost the same as the one from~\cite{Gu21a}. It looked like the estimate $3{.}5$ for the density in standard generators could not be exceeded. Many attepts to do that since 2004 led to an opinion that the construction by Belk and Brown was optimal. This was mentioned by Burillo in~\cite{Bur16}; the author discussed this with Jim Belk on some conferences. In~\cite{Gu22} we even formulated a conjecture that $i_*(F;A)=\frac12$ for the standard generating set $A=\{x_0,x_1\}$. However, it turned out that the conjecture was false.

In~\cite{Gu23} we got an improvement of the above estimate. The main result sounds as follows.

\begin{thm}[\cite{Gu23}]
	\label{over3.5}
	The density of the Cayley graph of Thompson's group $F$ in the standard set of generators $\{x_0,x_1\}$ is strictly greater than $3{.}5$. Equivalenly, the Cheeger isoperimetric constant of $F$ in the same set of generators is stricltly less than $\frac12$.
\end{thm}

The basic idea of the proof is as follows. The subset $BB(n,k)$ of the left Cayley graph of $F$ has some fragments with small density. The probability to meet such a fragment has a positive uniform lower bound. If we remove such fragments from the subgraph, we increase its density exceeding the value $3{.}5$. 

Let us describe the fragments we are interested in. Let we have a rooted binary forest $\dots,T_0,T_1,T_2,T_3,T_4,\dots$, where $T_0$ is marked and all trees have height at most $k$. We claim that trees $T_i$ ($0\le i\le4$) exist in this forest, and the following conditions hold:

\begin{itemize}
	\item $T_0$ and $T_2$ are trivial trees,
	\item $T_1$ and $T_3$ have height $k$,
	\item $T_4$ is a nontrivial tree.
\end{itemize}

The latter condition is added for simplicity. A forest satisfying the listed conditions is called {\em special}. To each of these forests we assign vertices $a$, $b$, $c$ of the left Cayley graph of $F$ in the standard generating set. Here $a$ corresponds to the forest with $T_0$ as a marked tree; $b$ and $c$ denote the forests where $T_1$ and $T_2$ are marked trees, respectively.

\begin{figure}[h]
\setlength{\unitlength}{.1cm}
\begin{center}
	\begin{picture}(87.75,18.552)(0,0)
		\put(23.75,13.75){\circle{9.014}}
		\put(46.25,13.75){\circle{9.014}}
		\put(69,13.5){\circle{9.014}}
		\put(42.25,13.75){\vector(-1,0){14.25}}
		\put(64.75,13.75){\vector(-1,0){14.25}}
		\put(19.5,13.75){\vector(-1,0){11.75}}
		\put(87.75,13.75){\vector(-1,0){14.5}}
		\put(46,9.5){\vector(0,-1){9.5}}
		\put(23.5,13.5){\makebox(0,0)[cc]{$a$}}
		\put(46.25,13.75){\makebox(0,0)[cc]{$b$}}
		\put(69.,13.25){\makebox(0,0)[cc]{$c$}}
		\put(13.75,16.75){\makebox(0,0)[cc]{$x_0$}}
		\put(36,16){\makebox(0,0)[cc]{$x_0$}}
		\put(57.75,15.5){\makebox(0,0)[cc]{$x_0$}}
		\put(80.75,15.5){\makebox(0,0)[cc]{$x_0$}}
		\put(50,3.25){\makebox(0,0)[cc]{$x_1$}}
	\end{picture}
	\vspace{1ex}
\end{center}

	\caption{The vertices $a,b,c$ in  the left Cayley graph of $F$ in standard generators.}\label{fig:picA}
\end{figure}
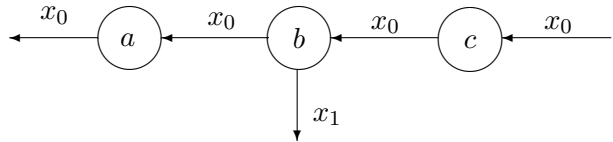

In the left Cayley graph of $F$ in standard generators these vertices look as in Figure~\ref{fig:picA}.

Since the tree $T_0$ is empty, we cannot remove a caret from it. This means that $a$ does not accept $x_1$ in the automaton $BB(n,k)$. The tree $T_1$ has height $k$ so a caret cannot be added to $T_0$ and $T_1$ to stay within $BB(n,k)$. So $a$ can accept only letters $x_0$ and $x_0^{-1}$. Exactly the same situation holds for the vertex $c$. Notice that the leftmost edge labelled by $x_0$ may not belong to the subgraph if $T_0$ is the leftmost tree in the special forest.

As for the vertex $b$, we can remove the caret from $T_1$ so $b$ accepts $x_1$. However, it does not accept $x_1^{-1}$ since the tree $T_1$ has height $k$ and no caret can be added to $T_1$ and $T_2$. Thus $a$, $c$ have degree 2 in $BB(n,k)$ and $b$ has degree 3 in the same subgraph.

If we have another special forest with the corresponding vertices $a'$, $b'$, $c'$, then no coincidences of vertices can occur. The only case could be $a'=c$ (or $a=c'$, which is totally symmetric). However, this is impossible by the choice of the tree $T_4$ in the special forest. If we go from $a'$ by a path labelled by $x_0^{-2}$, then we meet vertex $c'$ that corresponds to the trivial tree. Going along the path with the same label from $c$, we get the forest with marked tree $T_4$, which is nontrivial by definition. So this condition allows us to avoid repetitions.

Using properties of generating functions, we estimate the probability to meet a special forest. It turns out that the probability has a lower bound $p=\frac1{1200}$. This leads to finite subgraphs in the Cayley graph of $F$ with density exceeding $3{.}5004$. 

The size of the set we deal with is very huge: we estimate is as $2^{2^{7200}}$. Recall that, according to~\cite{Moore13}, the size of Folner sets in $F$ (provided it is amenable) has a very fast growth: as a tower of exponents. We think that Theorem~\ref{over3.5} increases the chances for the group $F$ to be amenable.

Our construction disproves one more conjecture from~\cite{Gu22}. We thought that the Cayley graph of $F$ for the set $\{x_0,x_1,x_2\}$ of generators has density $5$. In this case the Cheeger isoperimteric constant is $1$ so there exists a pure evacuation scheme on the graph according to Theorem~\ref{main2}. However, the following fact is true.

\begin{thm}[\cite{Gu23}]
	\label{dens012}
	The Cheeger isoperimetric constant of the Cayley graph of Thompson's group $F$ in the generating set $\{x_0,x_1,x_2\}$ is strictly less than $1$. Equivalently, the density of the corresponding Cayley graph strictly exceeds $5$.
\end{thm}

This means that there are no pure evacuation schemes on the Cayley graph of $F$ in these generators.
\vspace{1ex}

\begin{figure}[h]
\setlength{\unitlength}{.1cm}
\begin{center}
	\begin{picture}(87.75,18.553)(0,0)
		\put(23.75,13.75){\circle{9.014}}
		\put(46.25,13.75){\circle{9.014}}
		\put(69,13.5){\circle{9.014}}
		\put(42.25,13.75){\vector(-1,0){14.25}}
		\put(64.75,13.75){\vector(-1,0){14.25}}
		\put(19.5,14){\vector(-1,0){11.75}}
		\put(87.75,13.75){\vector(-1,0){14.5}}
		\put(46,9.5){\vector(0,-1){9.5}}
		\put(23.5,13.5){\makebox(0,0)[cc]{$a$}}
		\put(46.25,13.75){\makebox(0,0)[cc]{$b$}}
		\put(69.,13.25){\makebox(0,0)[cc]{$c$}}
		\put(13.75,16.75){\makebox(0,0)[cc]{$x_0$}}
		\put(36,16){\makebox(0,0)[cc]{$x_0$}}
		\put(57.75,15.5){\makebox(0,0)[cc]{$x_0$}}
		\put(80.75,15.5){\makebox(0,0)[cc]{$x_0$}}
		\put(50,3.75){\makebox(0,0)[cc]{$x_1$}}
		\put(23.25,9.75){\vector(0,-1){9.5}}
		\put(69,9.25){\vector(0,-1){9.5}}
		\put(27,3.75){\makebox(0,0)[cc]{$x_2$}}
		\put(72.25,3.75){\makebox(0,0)[cc]{$x_2$}}
	\end{picture}
	\vspace{1ex}
\end{center}

	\caption{Fragment of the left Cayley graph of $F$ in generators $\{x_0,x_1,x_2\}$ 
for any special forest.}\label{fig:pic}
\end{figure}
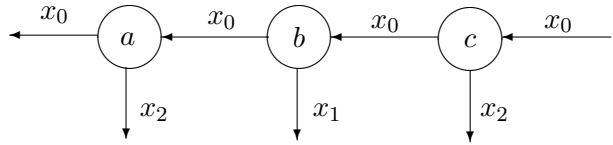

The fragment of the left Cayley graph of $F$ in generators $\{x_0,x_1,x_2\}$ will look as in Figure~\ref{fig:pic}
for any special forest (we keep previous notation).

Notice that trees $T_1$, $T_3$ of the special forest have height $k$. So no carets can be placed over any pair of trees of the form $T_i$, $T_{i+1}$, where $0\le i\le3$. This explains why no edges labelled by $x_1^{-1}$, $x_2^{-1}$ can be accepted by vertices in the picture. The vertex $a$ corresponds to a trivial marked tree so it does not accept $x_1$. However, it accepts $x_2$ since the tree to the right of $T_0$ has a caret. A similar situation holds for the vertex $c$. As for $b$, it accepts $x_1$ but it does not accept $x_2$ since the tree $T_2$ to the right of the marked tree $T_1$ is empty.

The condition on the tree $T_4$ allows us to avoid repetitions of vertices $a$, $b$, $c$ for different special forests.

Now we are subject to remove vertices of the form $a$, $b$, $c$ for all special forests. We also remove geometric edges incident to these vertices. For 3 vertices we remove no more than 14 directed edges from the graph. The key point here is inequality $\frac{14}3 < 5$. Here 5 is the limit of densities for $BB(n,k)$ considered as subgraphs in the left Cayley graph of $F$ in generators $\{x_0,x_1,x_2\}$. In the previous Section the same r\^ole was played by inequality $\frac{10}3 < 3{.}5$.

\section{Equations in group rings and their systems}
\label{eqgr}

Tamari~\cite{Ta54} shows that if a group $G$ is amenable, then the group ring $R=K[G]$ satsfies the Ore condition for
any field $K$. This means that for any $a,b\in R$ there exist $u,v\in R$ such that $au=bv$, where $u\ne0$ or $v\ne0$. 

This stament can be easily generalized. Suppose that instead of one linear equation $au=bv$ with coefficients in $R$ we have a system of them, where the number of variables exceeds the number of equations:
$$
\left\{
\begin{array}{ccccccc}
	a_{11}u_1&+  &\cdots  &+  &a_{1n}u_n  & = & 0 \\
	\cdots&  &  \cdots & &\cdots  &  &  \\
	a_{m1}u_1&+  &\cdots  &+  &a_{mn}u_n  &=  &0 
\end{array}
\right.
$$
where $n > m$, $a_{ij}\in R$ for all $1\le i\le m$, $1\le j\le n$. We are interested in solutions $(u_1,...,u_n)\in R^n$. 
\vspace{1ex}

Cardinality arguments based on the F\o{}lner criterion show that for amenable group $G$, this system always has a nonzero solution.
\vspace{1ex}

In~\cite{Bar} Bartholdi shows that the converse to the above statement is true. This gives a new criterion for amenabilty of groups. Although Theorem 1.1 in~\cite{Bar} concerns the so-called GOE and MEP properties of automata (Gardens of Eden and Mutually Erasable Patterns), the proof of it allows one to extract the following statement.

\begin{thm}[\cite{Bar}]
	\label{barth19}
	{\rm(Bartholdi)}\ For any group $G$, the following two properties are equaivalent.
	
	(i) $G$ is amenable
	
	(ii) For any field $K$ and for any system of $m$ linear equations over $R=K[G]$ in $n > m$ variables, there exists a nonzero solution.
\end{thm}

In the Appendix to the same paper, Kielak shows that if the group ring $K[G]$ has no zero divisors, both properties are equivalent to the Ore condition. In particular, this holds for R.\,Thompson's group $F$. It is orderable, so there are no zero divisors in a group ring over a field. So we quote the following

\begin{thm}[\cite{Bar}]
	\label{kielak}
	{\rm(Kielak)}\ The group $F$ is amenable if and only if the group ring $K[F]$ over any field satisfies the Ore condition.
\end{thm}

Let $M$ be a cancellative monoid. It is known from~\cite{Ta54} that if $M$ is left amenable then the monoid ring $K[M]$ satisfies Ore condition, that is, there exist nontrivial common right multiples for the elements of this ring. In~\cite{Don10} Donnelly shows that a partial converse to this statement is true. Namely, if the monoid $\mathbb Z^{+}[M]$ of all elements of $\mathbb Z[M]$ with positive coefficients has nonzero common right multiples, then $M$ is left amenable. He asks whether the converse is true for this particular statement.

In~\cite{Gu21b} we show that the converse is false even for the case of groups. Say, if $M$ is a free metabelian group, then $M$ is amenable but the Ore condition fails for $\mathbb Z^{+}[M]$. Besides, we study the case of the monoid $M$ of positive elements of Thompson's group $F$. Notice that the amenability problem for the group $F$ is equivalent to left amenability of the monoid $M$; on the other hand, it is known that $M$ is not right amenable~\cite{Gri90}. We show that for this case the monoid $\mathbb Z^{+}[M]$ does not satisfy Ore condition. That is, even if $F$ is amenable, this cannot be shown using the above sufficient condition.
\vspace{1ex}

The proof was based on the following

\begin{lem}[\cite{Gu21b}]
	\label{apm1bpm1}
	Let $M$ be a monoid embeddable into a group $G$. Suppose that the monoid $\mathbb Z^{+}[M]$ has nonzero common right multiples. Then for any $a,b\in M$ there exists a relation of the form $a^{\pm1}b^{\pm1}\ldots a^{\pm1}b^{\pm1}=1$ that holds in $G$.
\end{lem}

On the other hand, we prove that this property does not hold if $M$ is the free metabelian group with basis $\{a,b\}$. This implies

\begin{thm}[\cite{Gu21b}]
\label{thmdon1}
There exists a left amenable cancellative monoid $M$ (actually, an amenable group) such that the monoid $\mathbb Z^{+}[M]$ does not satisfy the Ore condition.
\end{thm}

The following property of the group $F$ clarifies the situation for it with respect to Donnelly's condition.

\begin{lem}[\cite{Gu21b}]
	\label{altrel}
	Let $x_0$, $x_1$, ... , $x_m$, ... be the standard generating set for R.\,Thompson's group $F$. Then any word of the form $w=x_{i_1}^{\pm1}x_{j_1}^{\pm1}\ldots x_{i_k}^{\pm1}x_{j_k}^{\pm1}$ $(k\ge1)$ does not represent the identity element of $F$ provided all $i_1$, ... , $i_k$ are even and all $j_1$, ... , $j_k$ are odd.
\end{lem}

Here is the consequence of the above statements:

\begin{thm}[\cite{Gu21b}]
\label{thmdon2}
Let $M$ be positive monoid of R.\,Thomspon's group $F$. Then $\mathbb Z^{+}[M]$ does not satisfy Ore condition.
\end{thm}

Now let us review the results from~\cite{Gu21c}. First of all, the group $F$ can be replaced by the monoid $M$ in the statement of the Ore condition:

\begin{lem}[\cite{Gu21c}]
\label{kfkm}
For any field $K$, the group ring $K[F]$ satisfies the Ore condition if and only if the monoid ring $K[M]$ satisfies the Ore condition.
\end{lem}

This means the we can just forget about negative powers of variables working in the ring of skew polynomials satisfting $x_jx_i=x_ix_{j+1}$ for $0\le i < j$. This follows from an elementary fact that any finite subset in $F$ can be sent into $M$ under right multiplication. More precisely:

\begin{lem}
\label{gig}
For any $g_1,...,g_n\in F$ there exists $g\in F$ such that $g_1g,\ldots,g_ng\in M$.
\end{lem}

Another step is to reduce the Ore equation $au=bv$ in $K[M]$ to the case when $a$, $b$ are homogeneous polynomials of the same degree.

\begin{lem}[\cite{Gu21c}]
\label{homogen}
Suppose that any equation of the form $au=bv$ has a nonzero solution in $K[M]$ provided $a$, $b$ are homogeneous polynomials of the same degree. Then $K[M]$ satisfies the Ore condition.
\end{lem}

Now we have a bunch of equations in $K[M]$ indexed by two parameters. The first one is $d$, the degree of homogeneous polynomials $a$ and $b$. The second one is $m$, where $m$ is the highest subscript in variables we involve. The general strategy can be the following: we try to solve as much equations in $K[M]$ as we can, using this classification. For a pair of numbers $d\ge1$, $m\ge1$, we can take $a$, $b$ as linear combinations of monomials of degree $d$ in variables $x_0$, $x_1$, ... , $x_m$ with indefinite coefficients. We can think about these coefficients as elements of the field of rational functions over $K$ with a number of variables.
\vspace{1ex}

More precisely, we can state the general problem as follows. Any finite system of monomials of degree $d\ge1$ is contained in a set of the form ${\mathcal S}_{m+1,m+d+1}$ for some $m\ge1$. This set consists of all elements in the monoid $M$ with normal forms $x_{i_1}\ldots x_{i_d}$, where $i_1\le\cdots\le i_d$ and $i_1\le m$, $i_2\le m+1$, ... , $i_d\le m+d-1$. Let $K[S]$ denote the set of all linear combinations of elements of $S\subset M$ with coefficients in $K$. (Strange parameters in the notation we use is explained by the following reason: if we represent elements of $M$ by positive semigroup diagrams, then the set of them will consists of $(x^{m+1},x^{m+d+1})$-diagrams.)
\vspace{1ex}

{\bf Problem} ${\mathcal P}_{d,m}$: {\sl Given two elements $a,b\in K[{\mathcal S}_{m+1,m+d+1}]$, find a nonzero solution of the equation $au=bv$, where $u,v\in K[M]$, or prove that it does not exist.}
\vspace{1ex}

According to Theorem~\ref{kielak} by Kielak, and Lemma~\ref{homogen} on homogeneous equations, we have the following alternative. If the Problem ${\mathcal P}_{d,m}$ has positive solution for any $d,m\ge1$ (that is, we can find nonzero solutions), then the group $F$ is amenable. If this Problem has negative solution for at least one case, then $F$ is not amenable.
\vspace{1ex}

We start with the case of polynomials of degree $d=1$ for arbitrary $m$. Let us consider equations of the form
\be{ab}
(\alpha_0x_0+\alpha_1x_1+\cdots+\alpha_mx_m)u=(\beta_0x_0+\beta_1x_1+\cdots+\beta_mx_m)v.
\ee
We show that all of them have nonzero solutions in $K[M]$ for arbitrary coefficients $\alpha_i,\beta_i\in K$ ($0\le i\le m$). This follows from cardinality reasons, see~\cite[Section 2]{Gu21c}. However, this approach gives us soultions of very high degree. More precisely, we prove the following

\begin{thm}[\cite{Gu21c}]
	\label{xmy}
	a) For any $m\ge1$, the set of elements $X_m=\{x_0,x_1,...,x_m\}$ is not doubling, that is, there exists a finite subset $Y\subset M$ such that $|X_mY| < 2|Y|$. 
	
	b) If $a,b\in K[M]$ are linear combinations of monomials $x_0$, $x_1$, ... , $x_m$ of degree $1$, then the equation $au=bv$ in $k[M]$ has a nonzero solution, where $\deg u=\deg v\le\frac{m(m+1)}2$.
\end{thm}

This means that Problem ${\mathcal P}_{1,m}$ has positive solution for any $m\ge1$. However, the degree of $u$, $v$ have quadratic growth with respect to $m$. This estimate can be essentially improved that will be shown later. Before doing that, we mention the following fact that follows from cardinality reasons.

\begin{thm}[\cite{Gu21c}]
	\label{s24}
	The equation of the form
	$$
	(\alpha_{00}x_0^2+\alpha_{01}x_0x_1+\alpha_{02}x_0x_2+\alpha_{11}x_1^2+\alpha_{12}x_1x_2)u=
	(\beta_{00}x_0^2+\beta_{01}x_0x_1+\beta_{02}x_0x_2+\beta_{12}x_1^2+\beta_{12}x_1x_2)v
	$$
	in the monoid ring $K[M]$ has a nonzero solution, where $\deg u,v\le41$.
\end{thm}

This means that Problem ${\mathcal P}_{2,1}$ has positive solution. The cases that come after that are already unknown. This is $d=2$, $m=2$, where $a$, $b$ are linear combinations of 9 monomials of degree 2:
$$
x_0^2,\, x_0x_1,\, x_0x_2,\, x_0x_3,\, x_1^2,\, x_1x_2,\, x_1x_3,\, x_2^2,\, x_2x_3.
$$
This is one possible candidate to obtain the negative answer. If, nevertheless, this Problem ${\mathcal P}_{2,2}$ has positive answer (that is, there exist nonzero solutions), then one can try Problem ${\mathcal P}_{3,1}$, where $a$, $b$ are linear combinations of 14 monomials of degree 3:
$$
x_0^3, x_0^2x_1, x_0^2x_2, x_0^2x_3, x_0x_1^2, x_0x_1x_2, x_0x_1x_3, x_0x_2^2, x_0x_2x_3, x_1^3, x_1^2x_2, x_1^2x_3,x_1x_2^2,x_1x_2x_3.
$$

The further strategy can be as follows: starting from equations of a simple form, we try not only to prove they have nonzero solutions (which is not so hard), but also try to describe somehow the set of all their solutions.

Say, if we have equation of the form $au=bv$, where $a=\alpha_0x_0+\alpha_1x_1$, $b=\beta_0x_0+\beta_1x_1$, then the description of all its solutions is easy. Namely, $u=(\beta_0x_0+\beta_1x_2)w$, $v=(\alpha_0x_0+\alpha_1x_2)w$ for any $w\in R=K[M]$. This means that the intersection of two principal right ideals $aR\cap bR$ is a principal right ideal.

We also know how to describe all solutions of the equation  $(\alpha_0x_0+\alpha_1x_1+\alpha_2x_2)u=(\beta_0x_0+\beta_1x_1+\beta_2x_2)v$. This will be done later. For this case, the description will be more complicated. In particular, the intersection $aR\cap bR$ is no longer a principal right ideal. 
\vspace{1ex}

If instead of one equation we have a system of equations of the form
$$
(\alpha_1x_0+\beta_1x_1)u_1=(\alpha_2x_0+\beta_2x_1)u_2=\cdots=(\alpha_kx_0+\beta_kx_1)u_k
$$
for any $k\ge2$, then it also has a nonzero solution. Indeed, the product $$(\alpha_1x_0+\beta_1x_1)(\alpha_1x_0+\beta_1x_2)\cdots(\alpha_kx_0+\beta_kx_{k+1})$$ 
is left divisible by  $\alpha_ix_0+\beta_ix_1$ for any $1\le i\le k$, which can be checked directly.
\vspace{1ex}

A much more interesting example of a system of equations looks as follows. Let us state it as a separate problem.
\vspace{1ex}

{\bf Problem} ${\mathcal Q}_k$: {\sl Given $k+1$ linear combinations of elements $x_0$, $x_1$, $x_2$, consider a system of $k$ equations with $k+1$ unknowns:
	$$
	(\alpha_0x_0+\beta_0x_1+\gamma_0x_2)u_0=(\alpha_1x_0+\beta_1x_1+\gamma_1x_2)u_1=\cdots=(\alpha_kx_0+\beta_kx_1+\gamma_kx_2)u_k.
	$$
Find a nonzero solution of this system, where $u_0,u_1,\ldots,u_k\in K[M]$, or prove that it does not exist.}
\vspace{1ex}

Notice that ${\mathcal Q}_1$ has been already considered. To solve ${\mathcal Q}_2$ in positive, it suffices to find a finite set $Y$ with the property $|AY| < \frac32|Y|$, where $A=\{x_0,x_1,x_2\}$. This can be done by cardinality reasons similar to the ones used in the proof of Theorem~\ref{s24}. The estimate there will be also $n\ge45$. In fact, we are able to prove a much stronger fact. Namely, using the result of~\cite{Gu21a}, we can construct a finite set $Y$ with the property $|AY| < \frac43|Y|$. The size of $Y$ is really huge, it does not have transparent description. This immediately implies that Problem ${\mathcal Q}_3$ has a positive solution. 
\vspace{1ex}

\begin{thm}[\cite{Gu22se}]
	\label{43}
	Let $A=\{x_0,x_1,x_2\}$. Then there exists a finite subset  $S\subset F$ satisfying $|AS| < \frac43|S|$.
\end{thm}

As an immediate corollary, we have

\begin{cor}[\cite{Gu22se}]
	\label{qq3}
	Let $R=K[F]$ be a group ring of $F$ over a field $K$. For any $4$ linear combinations of elements $x_0$, $x_1$, $x_2$ with coefficients in $K$, the system of $3$ equations with $4$ unknowns
	$$
	(\alpha_0x_0+\beta_0x_1+\gamma_0x_2)u_0=(\alpha_1x_0+\beta_1x_1+\gamma_1x_2)u_1=(\alpha_2x_0+\beta_2x_1+\gamma_2x_2)u_2=(\alpha_3x_0+\beta_3x_1+\gamma_3x_2)u_3
	$$
	has a non-zero solution in $R$.
\end{cor}

Donnelly shows in~\cite{Don14} that $F$ is non-amenable if and only if there exists $\varepsilon > 0$ such that for any finite set $Y\subset F$, one has $|AY|\ge(1+\varepsilon)|Y|$, where $A=\{x_0,x_1,x_2\}$ (see also~\cite{Don07}). For the set $Y$ here, one can assume without loss of generality that $Y$ is contained in ${\mathcal S}_{4,n}$ for some $n$. This gives some evidence that the amenability problem for $F$ has very close relationship with the family of Problems ${\mathcal Q}_k$. The case $k=4$ looks as a possible candidiate to a negative solution (that is, all solutions are zero). If true, this will imply that the constant $\varepsilon=\frac14$ fits into the above condition.
\vspace{1ex}

Now we give an improved version of Theorem~\ref{xmy}.

\begin{thm}[\cite{Gu21c}]
	\label{aubv}
	The equation
	$$
	(\alpha_0x_0+\alpha_1x_1+\cdots+\alpha_mx_m)u=(\beta_0x_0+\beta_1x_1+\cdots+\beta_mx_m)v
	$$
	has a nonzero solution, where $u$, $v$ are homogeneous polynomials of degree $m$ in variables $x_0,x_1,\ldots,x_{2m}$.
\end{thm} 

According to our strategy, we are interested in describing the set of all solutions to the above equation. This is easy to do for $m=1$. Let $m=2$. The equation from Theorem~\ref{aubv} can be rewritten as
\be{012}
(x_0+\alpha x_2)u=(x_1+\beta x_2)v,
\ee
up to a linear transformation, where $\alpha,\beta$ are some coefficients. We will assume both of them are nonzero: otherwise the description becomes trivial.

First of all, we take a solution of this equation extracted from the proof of Theorem~\ref{aubv}. One can check directly that the following polynomials satisfy~(\ref{012}):
$$
u_0=\beta x_0x_3+\beta^2x_0x_4-\alpha x_1x_3-\alpha\beta x_1x_4-\alpha\beta x_3^2-\alpha\beta^2 x_3x_4,
$$
$$
v_0=\beta x_0^2-\alpha x_0x_1-\alpha^2x_3^2-\alpha^2\beta x_3x_4.
$$
We say that $(u_0,v_0)$ is a {\em basic} solution. Now we are going to show how to extract all solutions from it.
\vspace{1ex}

By $M_1$ we denote the submonoid of $M$ generated by $x_1$, $x_2$, ... . We prove the following lemma using the same ideas as in the number-theoretical Remainder Theorem.

\begin{lem}[\cite{Gu21c}]
	\label{v0R}
	For any $v\in K[M]$ there exist $w_1\in K[M]$, $w_2,w_3\in K[M_1]$ such that $v=v_0w_1+x_0w_2+w_3$.
\end{lem}

Let us denote by $\phi$ an endomorphism of $F$ that takes each $x_i$ to $x_{i+1}$ ($i\ge0$). The description of the set of solutions for~(\ref{012}) can be reduced to the case when $\alpha=\beta$. Namely, if $(u,v)$ is a solution for~(\ref{012}) with $d=\deg u=\deg v > 1$, then the following equalities hold:
\be{eqs}
\left\{
\begin{array}{ccl}
	u	& = & u_0w_1+\alpha^{-1}(x_1+\beta x_3)\phi(u') \\
	v	& = & v_0w_1+\alpha^{-1}x_0\phi(v')+x_3\phi(u')
\end{array}
\right.
\ee
where $(u',v')$ is a solution of the equation~(\ref{012}) with $\alpha=\beta$, and $\deg u'=\deg v'=d-1$. 
So the problem can be reduced to the case when $\alpha=\beta$ in the equation. For this case we have the following inductive decription. 

\begin{thm}[\cite{Gu21c}]
	\label{descr}
	Let $\beta\ne0$ be an element of a field $K$. Let us consider the equation $(x_0+\beta x_2)u=(x_1+\beta x_2)v$ in the monoid ring $K[M]$. Let
	$$
	u_0=x_0x_3+\beta x_0x_4-x_1x_3-\beta x_1x_4-\beta x_3^2-\beta^2 x_3x_4,
	$$
	$$
	v_0=x_0^2-x_0x_1-\beta x_3^2-\beta^2x_3x_4
	$$
	be its basic solution.
	
	Then for any its solution, one has the following presentation for its first unknown:
	$$
	u=u_0w_0+(x_1+\beta x_3)\phi(u_0)w_1+(x_1+\beta x_3)(x_2+\beta x_4)\phi^2(u_0)w_2+\cdots+\prod\limits_{i=1}^k(x_i+\beta x_{i+2})\phi^k(u_0)w_k,
	$$
	where $k\ge0$, and $w_i$ belongs to $K[M_i]$, where $M_i$ is the submonoid of $M$ generated by $x_i,x_{i+1},\ldots\,$ $(0\le i\le k)$.
\end{thm}

Notice that if we know $u$, then the second unknown $v$ is determined uniquely since the group ring $K[F]$ over a field has no zero divisors.
\vspace{1ex}

Now let us review the results of~\cite{Gu22se} that continue the above research. We are going to consider equations and their systems in the group ring $K[F]$ of the group $F$ over a field $K$. Ring coefficients of the equations will not be assumed to be homogeneous polynomials in $K[M]$.
\vspace{1ex}

Let we have a group word $w=\xi_1\xi_2\ldots\xi_n$ where $\xi_i\in\{x_0^{\pm1},x_1^{\pm1}\}$. Suppose that $w=1$ in $F$. Let $g_i=\xi_i\ldots\xi_n$ for $1\le i\le n$ and $g_{n+1}=1$. We see that $(g_1-g_2)+(g_2-g_3)+\cdots+(g_n-g_{n+1})=g_1-g_{n+1}=w-1=0$ in the group ring $K[F]$. On the other hand, $g_i-g_{i+1}=(\xi_i-1)g_{i+1}$ for all $1\le i\le n$. Elements of the form $\xi_i-1$ are equal to $-(1-x_0)$, $-(1-x_1)$, $(1-x_0)x_0^{-1}$, $(1-x_1)x_1^{-1}$. Therefore, the above equation can be rewritten in the form $(1-x_0)u=(1-x_1)v$, where $u,v\in K[F]$. Elements $u$, $v$ are defined uniquely by the relation $w=1$ in $F$.

Recall that amenability of an $m$-generated group $G$ can be characterized in terms of its cogrowth rate. Let $A$ be a generating set for $G$, where $|A|=m$. Denote by $P_n$ the number of group words over $A$ of length $n$ representing the trivial element of $G$. Then $G$ is amenable if and only if $\limsup\limits_{n\to\infty}\sqrt[n]{P_n}=2m$. This is a famous Kesten criterion~\cite{Kest,Kest2}. If $\bar{P}_n$ the number of reduced group words over $A$ of length $n$ trivial in $G$, then $G$ is amenable if and only if $\limsup\limits_{n\to\infty}\sqrt[n]{\bar{P}_n}=2m-1$. This is Grigorchuk criterion~\cite{Gri80}.

According to these remarks, we see that it is useful to know the description of all solutions to the equation $(1-x_0)u=(1-x_1)v$ in the group ring of $F$. If we take one of the defining relations of $F$, namely, $x_1^{x_0^2}=x_1^{x_0x_1}$, and apply to it the above procedure, then we get to the following equation in the group ring of $F$:
$$
(1-x_0)\cdot(1-x_1)(1+x_1-x_2)=(1-x_1)\cdot(1-x_3-x_0^2+x_0x_1).
$$
We call the pair $(u,v)$ a {\em basic} solution of the equation $(1-x_0)u=(1-x_1)v$, where $u=(1-x_1)(1+x_1-x_2)$, $v=1-x_3-x_0^2+x_0x_1$. Now we are going to get all solutions in terms of the basic one. It suffices to present $u$ in its general form since $v$ is uniquely determined by $u$.

\begin{thm}[\cite{Gu22se}]
	\label{descr01}
	Let $R=K[F]$ be a group ring of $F$ over a field $K$. Let 
	$$
	u_0=(1+x_0-x_1)(1-x_3),
	$$
	$$
	u_1=(1-x_1)\phi(u_0)=(1-x_1)(1+x_1-x_2)(1-x_4),
	$$
	$$
	\dots
	$$
	$$
	u_k=(1-x_1)\phi(u_{k-1})=(1-x_1)\ldots(1-x_k)(1+x_k-x_{k+1})(1-x_{k+3})
	$$
	$$
	\dots
	$$
	
	Then for any solution of the equation $(1-x_0)u=(1-x_1)v$ in $R$, the element $u$ belongs to the right $R$-module generated by $u_0$, $u_1$, $u_2$, \dots\,. Moreover, for any solution in the monoid ring $K[M]$, one can express it as $u=u_0r_0+u_1r_1+\cdots+u_kr_k$ for some $k$ and elements $r_i\in K[M_i]$ $(0\le i\le k)$.
\end{thm}

It turns out that the equation $(1-x_0)u=bv$ has a non-zero solution for any element $b\in K[F]$.

\begin{thm}[\cite{Gu22se}]
	\label{x0b}
	Let $R=K[F]$ be a group ring of $F$ over a field $K$. Then for any element $b\in R$, the equation $(1-x_0)u=bv$ has a non-zero solution in $R$.
\end{thm}

The proof is based on the following fact that looks interesting in itself.

\begin{thm}[\cite{Gu22se}]
	\label{bphi}
	Let $b\in K[M]$ be any element in the monoid ring $R=K[M]$. Then the set $B=\{b,\phi(b),\phi^2(b),\cdots\}$ is not a free basis of the right $R$-module it generates.
\end{thm}

Notice that we do not know whether any equation of the form $(1-x_1)u=bv$ has a nontrivial solution. So we asked the following

\begin{qu}
	\label{x1b}
	Let $R=K[F]$ be a group ring of $F$ over a field $K$. Is it true that for any element $b\in R$, the equation $(1-x_1)u=bv$ has a non-zero solution in $R$?
\end{qu}

Using the automorphism $x_0\mapsto x_0^{-1}$, $x_1\mapsto\bar{x}_1=x_1x_0^{-1}$, one can reduce the above problem to the case of equation $(x_0-x_1)u=bv$, where $b\in K[M]$. 

We can extract from Theorem~\ref{x0b} the following

\begin{cor}[\cite{Gu22se}]
	\label{x0xm}
	For any $m\ge1$ there exists a non-zero solution to the system of equations
	\be{x0x1xm}
	(1-x_0)u_0=(1-x_1)u_1=\cdots=(1-x_m)u_m
	\ee
	in the group ring of $F$.
\end{cor}

The union of supports of the elements $u_0,u_1,\ldots,u_m\in K[F]$ is a multi-dimensional analog of a relation in $F$. An explicit form of this finite subset or even its size are unknown already for the case $m=2$.
\vspace{1ex}

One more general criterion for non-amenability of a group can be presented in terms of group series. Let $G$ be a group generated by a finite set $A=\{a_1,\ldots,a_m\}$. For any field $K$ we denote by $K[[G]]$ the space of infinite sums of the form
$$
\sum\limits_{g\in G}\alpha(g)\cdot g,
$$
where $\alpha(g)\in K$ are coefficients $(g\in G)$. Obviously, $G$ acts on the left and on the right on this space. These actions can be naturally extended to the group ring. This gives $K[[G]]$ the structure of a $K[G]$-bimodule. For the application we need, we assume that $K=\mathbb R$ will be the field of reals (or rationals, if necessary).
\vspace{1ex}

\begin{thm}[\cite{Gu22se}]
	\label{series}
	Let $G$ be a group generated by $A=\{a_1,\ldots,a_m\}$. Then $G$ is non-amenable if and only if there exists an equality in the left $\mathbb R[G]$-module $\mathbb R[[G]]$ of the form
	$$
	\sum\limits_{g\in G}g=(1-a_1)S_1+\cdots+(1-a_m)S_m,
	$$
	where $S_1,\ldots,S_m$ are elements of $\mathbb R[[G]]$ with uniformly bounded coefficients.
\end{thm}

In partucular, this is applied to $F$, where we need an equality of the form $\sum\limits_{g\in F}g=(1-x_0)u+(1-x_1)v$, where $u$, $v$ are infinite group series from $\mathbb R[[F]]$ with uniformly bounded coefficients.

\section{Miscellaneous}
\label{misc}

In this section we describe several results from~\cite{Gu22s} concerning the group $F$. They involve a simple description of exhausting finite fragments of the Cayley graphs as well as a new algorithm for solving the word problem in this group.
We have already mentioned that, in the context of the amenability problem, instead of describing the infinite Cayley graph it is sufficient to describe finite subgraphs that contain balls of arbitrarily large radius. Moreover, it is enough
to consider such finite fragments within the monoid $M$ of positive elements.

\subsection{Structure of Cayley graphs}
\label{scg}

One of the inductive descriptions of finite fragments of the Cayley graph was given in~\cite{Gu04} in connection with producing a family of subgraphs whose density approaches $3$. A fairly standard tool for working with the group is also a representation of its elements as marked rooted binary forests, which we used in Section~\ref{dencg} in the description of
the Belk-Brown construction. First of all, we will show how to associate with every marked rooted binary forest an ordered triple of trees, which can be done in a fairly natural way.

Suppose that we are given a marked rooted binary forest, the trees of which are
denoted as $T_{-k}$, \dots, $T_{-1}$, $T_0$, $T_1$, \dots, $T_m$, where $k,m\ge0$ and the marked tree is $T_0$.

We add to this forest one-point trees on the left and on the right, denoting them by $T_{-(k+1)}$ and $T_{m+1}$, respectively. We further add carets, attaching them consecutively to the two rightmost trees without touching the tree $T_0$ and proceeding until there is exactly one tree on the right of $T_0$. In doing so we add $m$ carets.
We conduct exactly the same procedure on the left, adding $k$ carets. Namely, we consecutively attach a new caret to the two leftmost trees until there is exactly one tree on the left of $T_0$. The result of this process is a triple of rooted binary trees. The tree that was initially marked remains in the middle in an unaltered form. 

Let $n$ denote the total number of leaves in the triple of trees. Clearly, $n\ge3$, and it is easy to see from elementary combinatorial considerations that the number of triples of trees for this $n$ can be expressed in terms of the Catalan numbers as
$$
c_{n-1}-c_{n-2} = \frac{3(2n-4)!}{n!(n-3)!}.
$$

Let us define a graph $\Gamma_n$. Its vertices are the ordered triples of rooted binary trees described above with total number of leaves equal to $n\ge3$. The edges correspond to the actions of the generators that we shall now describe.

Multiplication of an element of the group on the left by $x_0$ corresponds to the following operation: we remove the leftmost caret, obtaining a quadruple of trees, and place a caret on the two rightmost trees. Such an operation is applicable if the first tree of the triple is non-trivial. Multiplication on the left by $x_0^{-1}$ is defined
similarly. Multiplication of an element of the group on the left by $x_1$ corresponds to the following operation: we remove the middle caret, obtaining a quadruple of trees, and place a caret on the two rightmost trees. Such an operation is applicable if the second tree of the triple is non-trivial. Multiplication on the left by $x_1^{-1}$ is defined
similarly.
 
Here we can observe that the operation of multiplication on the left by the symmetric generator $\bar{x}_1=x_1x_0^{-1}$ corresponds to removing the middle caret and then moving it to the left rather than to the right. In this example we clearly see this type of symmetry. It corresponds to the outer automorphism of $F$ defined by the rule $x_0\mapsto x_0^{-1}$, $x_1\mapsto\bar{x}_1$. 

Thus, $\Gamma_n$ becomes a subgraph of the left Cayley graph of $F$. The procedure of moving the carets (the application of the generators in the set $\{x_0^{\pm1},x_1^{\pm1},\bar{x}_1^{\pm1}\}$) somewhat resembles the well-known “Tower of Hanoi”.

The structure of these graphs is quite convenient for working with the group, especially for constructing flows on the Cayley graphs. As we already know, it is sufficient to do this on expanding finite fragments. The verification of the fact
that graphs of the form $\Gamma_n$ contain balls of arbitrarily large radii is fairly easy.

We proceed to describe a simpler model for fragments of the Cayley graph. In the graph $\Gamma_n$ the vertices are fairly complicated geometric objects, namely, triples of trees. We shall now achieve a situation where the vertices of the graph are ordinary geometric points. The edges will have an equally simple look. We consider the disjoint union of all rooted binary trees with a given number of carets equal to $n$. Then we have $c_n$ such trees. Consider an arbitrary vertex $v$ of
one of the trees. If it is not the uppermost one, then it occurs as part of exactly one caret situated over this vertex. We paint this caret red. If its top point is not the root of the tree, then we again consider the caret situated above this vertex and also paint it red. We continue the process until we reach the root of the tree. We begin removing the painted set of carets proceeding downwards from the top until we reach $v$. We shall obtain a rooted binary forest, which we turn into a marked forest by distinguishing the tree with root $v$. Such a procedure can be followed for
any vertex of any of the trees. 

We say that two vertices are equivalent if the procedure described above results in the same picture, that is, in the same marked rooted binary forest. Obviously, this relation $\sim$ on the set of vertices of all the trees under consideration is indeed an equivalence. (It is easy to observe that the total number of vertices in our disjoint union is equal to $(2n+1)c_n$.) We take all the available vertices and identify equivalent ones, obtaining the vertices of a new graph. The edges will have the following structure: for every caret consisting of a left and a right segment, we place arrows on these line segments proceeding downwards from the top and mark the left segment with the letter $x_1$ and the right segment with $\bar{x}_1$. It is understood that the edges connecting pairs of respectively equivalent vertices are also identified. What we have constructed is a subgraph of the left Cayley graph of $F$ in the symmetric system of group generators $\{x_1,\bar{x}_1\}$ in a, so to speak, ‘scattered’ form. The same vertex of the graph can occur as many copies, and the same applies for the edges. We can trace how many times the same vertex $v$ appears in the trees under
consideration. If, after removing the upper carets, we obtained a forest in which $s$ trees are situated to the left of $v$, and there are $t$ trees to the right of $v$, then the vertex $v$ occurs as $\frac{(s+t)!}{s!t!}$ copies. This fact is easily proved by induction on $s+t$ by using the properties of binomial coefficients.

We summarize what has been said above in the form of a separate assertion.

\begin{thm}[\cite{Gu22s}]
\label{easygraph}
The disjoint union of all trees with $n\ge0$ carets and with the equivalence relation $\sim$ on the set of its vertices defines a subgraph of the left Cayley graph of the group $F$ with the set of generators $\{x_1,\bar{x}_1\}$, which is isomorphic to the graph $\Gamma_{n+3}$. Furthermore, all the left segments of the carets have label $x_1$, and the
right ones have label $\bar{x}_1$ (in the direction from the top of the caret).

\end{thm}

We point out that edges with label $x_0$ do not, in fact, disappear anywhere: they can be added if we associate a triangle with every caret, and then an edge marked $x_0$ will go from right to left. We also emphasize that paths in the Cayley graph can now be depicted on the trees themselves, and we can ‘jump’ from any vertex to an equivalent one to
continue the path. Earlier, for an image of a path in the graph $\Gamma_n$, we had to draw all the vertices occurring in it in the form of complicated pictures, whereas this is no longer necessary in the new model.
\vspace{1ex}

\subsection{Algorithm for the word problem}
\label{awp}

We now proceed to describing a new algorithm for solving the word problem in $F$ without going beyond the group alphabet $\{x_0,x_1\}$ and without changing the length of the word being analyzed. We state everything in the form of a separate assertion. Recall that if a word is equal to $1$ in $F$, then it belongs to the derived subgroup of the free group (the sum of exponents both in $x_0$ and in $x_1$ is equal to zero). This follows from the fact that the defining relations have this property. Therefore we consider only such words.

\begin{thm}[\cite{Gu22s}]
\label{newalg}
Let $w$ be a group word over the alphabet $\{x_0,x_1\}$ with zero sum of the exponents of each of the generators. We describe a transformation ${\mathcal T}$ of $w$ at every step. 

We write the word in a circle representing it in the form of a labelled cyclic graph. We go around this graph in one of the directions starting from an arbitrary vertex, to which we assign weight $0$. Each consecutive vertex has the same weight as the preceding one if we go over an edge with label $x_1^{\pm1}$. When we go over an edge with label $x_0$, the weight of the next vertex increases by $1$. When we go over an edge with label $x_0^{-1}$, the weight of the next vertex decreases by $1$.

If a word does not contain the letters $x_1^{\pm1}$, then it is equal to $1$ in the free group, and therefore it is equal to $1$ in $F$. Otherwise, we consider edges with label $x_1^{\pm1}$. Their endpoints have equal weight. We choose those edges with label $x_1^{\pm1}$ for which this weight is maximal. We replace their labels according to the rule $x_1\mapsto x_0$, $x_1^{-1}\mapsto x_0^{-1}$. This yields a word of the same length, which is denoted by ${\mathcal T}(w)$.
	
If the sums of the exponents of $x_0$ and $x_1$ are zero in the new word, then we iterate the process by applying the transformation ${\mathcal T}$ again. If this is not the case, then $w$ is not equal to $1$ in $F$. The process terminates after finitely many steps. For a word equal to $1$ in $F$, at the end we must obtain a word only in $x_0^{\pm1}$.
\end{thm}

We can give a simple illustration of how the algorithm in Theorem~\ref{newalg} works. Consider the word $w=x_0^{-2}x_1x_0^2x_1^{-1}x_0^{-1}x_1^{-1}x_0x_1$ (one of the defining relations). It is easy to verify that the second and fourth occurrences of $x_1^{-1}$ have the maximum weight. After they are renamed as $x_0^{\pm1}$, the balance of exponents is preserved. The word ${\mathcal T}(w)$ takes the form $x_0^{-2}x_1x_0^2x_0^{-1}x_0^{-1}x_1^{-1}x_0x_0$ (we can cancel in the free group, but do not have to), and the two remaining occurrences of $x_1^{\pm1}$ acquire the
same (and therefore maximum) weight. At the next step they are renamed as $x_0^{\pm1}$.

One can hope that this algorithm will somehow help to estimate the number of words of given length that are equal to $1$ in $F$. In particular, it is possible to estimate the number of words that require a given number of steps of the algorithm. If a successful estimate could be obtained, using the Kesten -- Grigorchuk criterion would solve the amenability problem either in one way or the other.
	
\subsection{Perspectives}
\label{persp}

Here we discuss possible ways to work with the amenabilty problem for $F$. There can be three different strategies: proving non-amenability; proving amenabilty; investigation	in both directions.
\vspace{1ex}

1) One way to establish non-amenability is to find an evacuation scheme on the Cayley graph of $F$. This can be done inductively using the exhausting finite fragments desribed in Section~\ref{scg}. Say, this can be a sequence of subgraphs $\Gamma_n$ whose vertices are triples of trees. From Theorem~\ref{isopsym0} we know that there are no pure evacuation schemes for corresponding generating sets. Since it is easier to work with pure evacuation schemes rather than arbitrary ones, it is natural to replace the sets of triples by the set of, say, quadruples of trees in order to get the generating set wider.
\vspace{1ex}

2) Another way to establish non-amenability is to find an equation $au=bv$ in the group ring of $F$ that has only zero solutions. Possible candidates for that were given in Section~\ref{eqgr}. In this case we usually have to know a full description to a set of solutions for some classes of simpler equations. This looks like Fermat's proof that the equation $X^4+Y^4=Z^2$ has no solutions in positive integers based on the description of Pythagorean triples.
\vspace{1ex}

3) One more way is suggested by a paper by Wajnryb and Witowicz (published in the arXiv several years ago and then withdrawn by the authors). This was quite  an original approach. Let $G$ be a group generated by a finite set $A$. Consider any cyclic order on the set $A^{\pm1}$. A word of the form $xy^{-1}$ of length $2$ is called {\em regular} whenever $x$ is followed by $y$ in the cyclic order. The authors proved a general property: if $G$ is amenable, then for any $\gamma\in(0,1)$ there exists a relation of length $n$ in the group $G$ such that the number of regular subwords of consecutive letters is at least $\gamma n$.

The authors tried to prove that this property fails for $F$ for the cyclic order of letters $x_0$, $x_1^{-1}$, $x_0^{-1}$, $x_1$. The proof was incorrect, but the question whether $F$ satisfies the above property remains open. One can try any order for any generating set of $F$. If one can show that any relation of length $n$ in $F$ between the generators contain no more than $\gamma_0n$ regular subwords for some $\gamma_0 < 1$, this would imply that $F$ is not amenable.
\vspace{1ex}

4) To establish amenability, it suffices to construct a sequence of Folner sets in $F$. Roughly speaking, these are finite sets with high density. According to our Theorem~\ref{over3.5}, we know that the estimate $3{.}5$ for the density is not the best one. This increases the chances for $F$ to be amenable. So one needs to guess how Folner sets may look like. Some information can be found in~\cite{Moore13}. We also mention our recent paper~\cite{Gu23b}, where we consider a natural partition of $F$ into $7$ subsets according to the structure of normal forms (or canonical diagrams). We show that all but one of these sets may have non-zero measure provided finitely additive right invariant probability measure exists on $F$. Elements of zero-measured sets can be excluded from any Folner sets. This can help to understand the rule of constructing the desired sets. Roughly speaking, we need to know what elements of the group we want to take into the sets and what elements will be out of our sets.
\vspace{1ex}

5) One more way to prove amenability is to show that any equation $au=bv$ in the group ring of $F$ has a non-zero solution. This is equivalent to amenability according to Theorem~\ref{kielak} by Kielak. In this case we may try to solve Problem ${\mathcal P}_{d,m}$ for any parameters; first steps in this direction were considered in Section~\ref{eqgr}. Equally, one can consider any systems of equations in the monoid ring $K[M]$, where the number of unknowns exceeds the number of equations. In this case an inductive strategy becomes possible. Say, for any equation $au=bv$ one can decompose the elements by powers of $x_0$ with coefficients from $K[M_1]$. This leads to a system of equations in $K[M_1]$ with simpler coefficients.
\vspace{1ex}

6) A universal approach, when we do not assume the answer (positive or negative) can be based on Kesten -- Grigorchuk criterion. Here we need to estimate the number of group relations in $F$ of given length between $x_0$ and $x_1$. One tool can be extracted from Theorem~\ref{descr01}, where we describe all solutions to the equation $(1-x_0)u=(1-x_1)v$ in the group (or monoid) ring. We already mentioned that to any relation between the generators one can naturally assign such a solution.

From the same point of view, we have to mention Section~\ref{awp}, where we give an algorithm to solve the word problem in $F$ standing inside the set of group words over the alphabet $\{x_0,x_1\}$. Standard algorithms to solve the word problem based on the normal forms of elements involve an infinite alphabet $\{x_0,x_1,x_2,\ldots\}$, where the situation looks more difficult to control. Using Theorem~\ref{newalg}, we have a hope to get a lower or an upper bound for the number of relations in $F$ of given length.

Notice that even Theorem~\ref{easygraph} gives some light to the description of relations between symmetric generators $x_1$, $\bar{x}_1$.

\end{document}